\documentclass[11pt]{article}
\usepackage{setspace}
\usepackage[colorlinks=true]{hyperref}
\usepackage{ams math}
\usepackage{amssymb}
\title{Nov.25,2020  Quadratic Dynamics Over Hyperbolic Numbers}

\vspace{0.5 cm}
\begin{document}

 \hspace{0.5cm}  { QUADRATIC  DYNAMICS OVER HYPERBOLIC NUMBERS}\\                    
                           
                            \vspace{0.5 cm} 
                       \hspace {3.5cm} Sandra Hayes

      \vspace{0.5 cm} 
      
      \hspace{3.7
      cm} {\footnotesize $\mathbf {Abstract }$}\\     
          
          \vspace{0.5cm} 
                              
 {\footnotesize Hyperbolic numbers are a variation of complex numbers, but their dynamics is quite different. The hyperbolic Mandelbrot set for quadratic functions over hyperbolic numbers is simply a filled square, and the filled Julia set for hyperbolic parameters inside the hyperbolic Mandelbrot set is a filled rectangle. For hyperbolic parameters outside the hyperbolic Mandelbrot set, the filled Julia set has 3 possible topological descriptions, if it is not empty, in contrast to the complex case where it is always a non-empty totally disconnected set. These results were proved in [1,2,4,5,6,7] and are reviewed here.}
      
       \vspace{1.0cm}

The boundary of the Mandelbrot set for quadratic functions over complex numbers is a classical fractal. However, the Mandelbrot set for quadratic functions over hyperbolic numbers is simply a filled square, The filled Julia set for all hyperbolic parameters inside the hyperbolic Mandelbrot set is a filled rectangle, but the only filled Julia set for complex numbers whose boundary is a simple curve is when the parameter is the origin. For hyperbolic parameters outside the hyperbolic Mandelbrot set, the filled Julia set has 3 possible topological descriptions in contrast to the complex case where the filled Julia set for parameters outside the Mandelbrot set is always a Cantor set. \\     
       These known results, whose proofs are straightforward using the technique of characteristic coordinates, will be discussed here. In summary, the quadratic dynamics over hyperbolic numbers is uninteresting, in contrast to the importance of hyperbolic numbers in physics where they are viewed as basic to Einstein's theory of special relativity [9]. Some open questions for further investigations of the dynamics over hyperbolic numbers will be mentioned. \\        
     
      Definition:  A $ \emph{hyperbolic\ number}$ is a number $z$ with two real components $x$ and $y$ written as $ z = x+jy $ with a new imaginary unit $j$, called the $\emph{hyperbolic \ imaginary\ unit}$,  defined by the property that it is not equal to $1$ or $-1$ but satisfies $j^2 = 1$. Let $\mathbf{H}$ denote the set of all such numbers.\\
     
     Whereas the complex numbers form a field, the hyperbolic numbers are only an algebra, because there are zero divisors, i.e. non-zero hyperbolic numbers which have no multiplicative inverse. The zero divisors are all non-zero numbers $z = x + j y$ with $x = \pm y$, i.e  numbers on either diagonal. The zero divisors play a fundamental role in explaining the dynamics of quadratic functions over hyperbolic numbers as well as in explaining the applications  of hyperbolic numbers to physics. Geometrically, the zero divisors form the light cone. \\
    
     Because of the zero divisors, the algebra $\mathbf{H}$ is not a normed space in the natural way, i.e. by generalizing the Euclidean norm for complex numbers $||z|| = \sqrt{z \overline z} = \sqrt{x^2 + y^2}$, where $ \overline z = x-iy$ is the complex conjugate of the complex number $z = x + i y$. However, the $\emph{ modulus} $ of a hyperbolic number $z = x + jy $, which is defined to be $ \sqrt{|x^2 - y^2|}$, does play a role in special relativity, since it is invariant under the Lorentz transformations of special relativity.\\

      The hyperbolic Mandelbrot set as well as the hyperbolic filled Julia set for  quadratic functions $f_c(z) = z^2 + c$ over hyperbolic numbers, i.e.  $z$ and $c$ are hyperbolic, rely on the notion of boundedness for a sequence $(z_n)_{n \in \mathbf{N}}$ of hyperbolic numbers $z_n$.\\
   
    Definition:   $(z_n)_{n \in \mathbf{N}}$ is $\emph{bounded}$ if each of the real sequences $(x_n)_{n \in \mathbf{N}}$ and $(y_n)_{n \in \mathbf{N}}$ is bounded  where $z_n = x_n + jy_n$.\\

  This is one way of generalising the notion of boundedness for a sequence  $(z_n)_{n \in \mathbf{N}}$ of complex numbers  $z_n = x_n + iy_n$ to that for a sequence of hyperbolic numbers, but is not the only way as will be mentioned later.\\

 Recall that the complex Mandelbrot set  $\mathbf{M}$ comprises all complex numbers $c$ for which the orbit $(f_c^n(0))_{n \in \mathbf{N}}$ is bounded when $f_c(z) = z^2 + c$ with $z$ complex. Analogously, the $emph{hyperbolic\ Mandelbrot\ set }$, denoted by  $\mathbf{M_H}$, is the set of all hyperbolic numbers $c$ such that the forward orbit of the origin $(f_c^n(0))_{n \in \mathbf{N}}$ is bounded when $f_c(z) = z^2 +c$ and $z$ is hyperbolic.. \\
 
The technique of characteristic coordinates is used to prove that  $\mathbf{M_H}$ is a square. This is a tool using the existence of zero divisors and therefore not available in studying complex numbers.  
Define the $\emph{hyperbolic\ conjugate}$ of a hyperbolic number $z =x+jy$ by  $z ^*= x  -jy $. The zero divisors are all numbers $z = x + jy$ with $zz^* = x^2 - y^2 = 0$, and the modulus of $z$ is $ \sqrt{|z z^*|}$.
\\

 The $\emph{ characteristic\ coordinates}$ for the hyperbolic number $z=x+jy$ are  the real numbers $X = x - y, Y = x + y$. These are the  coordinates for $z$ with respect to the $\emph{ idempotent \ basis }$  $\{\alpha, \alpha^*\}$    
 with  $\alpha= \frac{1}{2}( 1- j)$ and $\alpha^* = \frac{1}{2}(1+j) $. Note that $\alpha$ and $\alpha^*$  are zero divisors. In this basis,  $ z = (x -y) \alpha + (x +y) \alpha^* = X \alpha + Y \alpha^*$, and the quadratic form $zz^* = x^2 - y^2$ for $z$ becomes simply the product of its characteristic coordinates $ XY,$ i.e. $zz^* = XY$. The term "idempotent" is justified, since  $\alpha^2 = \alpha, (\alpha^*)^2 = \alpha^*$. \\

 Let $f_c(z) = z^2 + c$ with hyperbolic numbers $z = x+ jy$ and $c = a + jb$. Then  $f_c(z) = x^2 + y^2 + a + j(2 xy + b)$ and its characteristic coordinates are $X^2 + c_1$ and $Y^2 + c_2$  where $c_1 = a-b$   and $c_2 = a + b$ are the characteristic coordinates of $c$. \\
 
 Consider the isomorphic linear transformation $T$ of  the vector space $\mathbf{H}$     
 given by the matrix\\
 \begin{equation}
 T =
  \begin{pmatrix}
 1& -1\\
  1 & 1
  \end{pmatrix}\\
  \end{equation}

which maps a point onto its characteristic coordinates i.e. $T(z)  = X + jY$, when hyperbolic numbers $z = x + jy$ are identified with ordered pairs $z = \begin{pmatrix} x\\ y \end{pmatrix}$ of real numbers.\\
Then \\

$T(f_c(z))  = X^2 + c_1+j( Y^2 +c_2) = f_{c_{1}} (X) +j(f_{c_{2}}(Y))$\\

for the real quadratic functions  $f_{c_1} (X) = X^2 + c_1,
 f_{c_2} (Y) = Y^2 + c_2$. By induction,  for all $n$\\

$T(f_c^n(z)) = f_{c_1}^n(X)+j( f_{c_2}^n(Y))$. \\

Therefore, orbits $f_c^n(0)$ of the origin for quadratic functions of hyperbolic numbers $c$ are reverted back to orbits $f_{c_1}^n(0)$ and $ f_{c_2}^n(0))$ of the origin for real quadratic functions which have been extensively studied. Notice that a hyperbolic sequence $z _n = x_n + jy_n$ is bounded if and only if its image $T(z_n)$ is bounded.\\
 
 The hyperbolic Mandelbrot set is a square as was shown in [1,2,4,5,6]: \\

 THEOREM 1.  For the quadratic hyperbolic function $f_c(z) = z^2 + c$ with  $c = a+jb $, let  $S= \{ (a, b) \in \mathbf{R^2} : -2 \leq a-b \leq \frac{1}{4}, -2 \leq a+b\leq \frac{1}{4} \} .$ Then\\
 
 \hspace{ 3.0cm} $\mathbf{M_H} = S $.\\

  PROOF: By definition, $c \in \mathbf{M_H}$ if and only if the sequence of hyperbolic numbers $(f_c^n(0))_{n \in \mathbf{N}}$ is bounded, i.e. if and only if both of the real sequences  $(f_{c_1}^n(0))_{n \in \mathbf{N}}$ and $(f_{c_2}^n(0))_{n \in \mathbf{N}}$ are bounded when $c_1 = a-b, c_2 = a + b$ are the characteristic coordinates of $c$. That means that $c_1$ and $c_2$ are in the real part of the complex Mandelbrot set: \\
 
  \hspace{2cm}   $\mathbf{M}\cap \mathbf{R} = [-2, \frac{1}{4}]$. \\
    
   Thus,  $ c_1, c_2  \in [-2, \frac{1}{4}] $ which defines $S$. \\
   
   Visually, in the ab-plane, one diagonal of the square $S$ is $[-2, \frac{1}{4}]$ on the a-axis, its  side length is $\frac{9}{8}\sqrt{2}$, and it intersects the b-axis at $\pm\frac{1}{4}$.\\
   
   Another definition of boundedness was used in [2], namely a sequence  $(z_n)_{n \in \mathbf{N}}$ of hyperbolic numbers $z_n = x_n + jy_n$ is $\Large{bounded}$ if $(|z_nz_n^*|)_{n \in \mathbf{N}}= ( |x_n^2 - y_n^2|)_{n \in \mathbf{N}}$ is bounded. That implies that every sequence of zero divisors  is bounded. In [2], the hyperbolic Mandelbrot set is, therefore, larger, namely\\
  
   \hspace{ 2.0cm} $\mathbf{M_H} = S \cup D$,\\
   
   where $D$ represents the two diagonals. \\

   \vspace{1.0 cm} 
        The complex filled Julia set $ \mathbf{K} (f_c) $ for the complex quadratic function $f_c(z) = z^2 + c $ is the set of all complex numbers $z$ for which $(f_c^n(z))_{n \in \mathbf{N}}$ is bounded. This set is either connected or a Cantor set and thus totally disconnected, but it is never empty since it contains all periodic points. \\
         Similarly, the  $\emph{hyperbolic\ filled\ Julia\ set }$ $\mathbf{K}_{\mathbf{H}}(f_c)$ for the hyperbolic quadratic function  $f_c(z) = z^2 + c$  is the set of all hyperbolic numbers $z = x + j y$ whose forward orbit $(f_c^n(z))_{n \in \mathbf{N}}$ is bounded or equivalently such that its image $(T(f_c^n(z))_{n \in \mathbf{N}}  $ is bounded, i.e  such that the sequences $(f_{c_1}^n(X))_{n \in \mathbf{N}}$ and $( f_{c_2}^n(Y))_{n \in \mathbf{N}}$ of real quadratic functions are bounded where $X = x-y, Y = x+y$ are the characteristic coordinates of $z$ and $c_1 = a-b, c_2 = a+b$  are the characteristic coordinates of $c = a + j b$. This set can be connected, totally disconnected,  disconnected but not totally disconnected or empty. The following two lemmas are used in the proof. \\  
     
     If $c$ is real, let \\

 \hspace{2.5cm}  $ \mathbf{K}_{\mathbf{R}}(f_c) = \mathbf{K} (f_c)\cap \mathbf{R}$\\

 be the set of all real numbers $x$ whose orbit $(f_c^n(x))_{n \in \mathbf{N}}$ is bounded. The hyperbolic filled Julia set can be written in terms of these real filled Julia sets [2, 6]: \\

      LEMMA 1.  Let $c = a+ j b $ with its characteristic coordinates $c_1 = a - b, c_2 = a + b $. If $f_c(z) = z^2 + c$ is a hyperbolic quadratic function, let  $f_{c_1}(X) = X^2 + c_1$ and $f_{c_2}(Y) = Y^2 + c_2 $ be the corresponding real quadratic functions in the characteristic coordinates $X = x-y, Y= x+y$ for $z = x + j y$. Then \\
       
  \hspace{2.0cm} $T( \mathbf{K}_{\mathbf{H}}(f_c)) = \mathbf{K}_{\mathbf{R}}(f_{c_1}) +j( \mathbf{K}_{\mathbf{R}}(f_{c_2}))$.\\

       PROOF:   $ z  \in  \mathbf{K}_{\mathbf{H}}(f_c)$ if and only if  the sequence $(f_c^n(z))_{n \in \mathbf{N}}$ is bounded which is true if and only if the real sequences        
     $ (f_{c_1}^n(X))_{n \in \mathbf{N}} $ and  $(f_{c_2}^n(Y))_{n \in \mathbf{N}}$ are bounded.  It follows that  $ z  \in  \mathbf{K}_{\mathbf{H}}(f_c)$ if and only if  $X \in \mathbf{K}_{\mathbf{R}}(f_{c_1})$ and  $Y\in\mathbf{K}_{\mathbf{R}}(f_{c_2} )$.\\
       
    To describe the hyperbolic filled Julia sets topologically,  the following well known facts [3] about the real filled Julia set of the real quadratic function $f_c(x) = x^2 + c$ for $x,c \in\mathbf{R}$ are used, namely it is empty or it is an interval or it is a Cantor set and therefore a totally disconnected subset of that interval.\\

Let $p_{\pm} = \frac{1 \pm \sqrt{1-4c}}{2}$ denote the two fixed points of the complex function $f_c$ which are real if and only if $c \leq \frac{1}{4}$. If $c = \frac{1}{4}$ there is just one fixed point at $p_{\pm} = \frac{1}{2}$.\\
          
      LEMMA 2.  Let $f_c(x) = x^2 + c$ for $x,c \in\mathbf{R}$.\\
       If  $c<-2$, then  $\mathbf{K}_{\mathbf{R}} (f_c) $ is a Cantor set  contained in $[-p_+, p_+].$ \\
             If  $c\in [-2, \frac{1}{4}]$, then $\mathbf{K}_{\mathbf{R}} (f_c) = [-p_+, p_+]$. \\
       If  $c>\frac{1}{4}$, then $\mathbf{K}_{\mathbf{R}} (f_c) = \emptyset$.\\
       
After Lemma 1, the hyperbolic filled Julia set is given by an ordered pair of two real filled Julia sets, each of which can be connected or a Cantor set, and thus totally disconnected, if it isn't empty. Then there are only three ways of combining them if neither is empty: both connected, both totally disconnected or one connected and the other totally disconnected. The following topological description of the hyperbolic filled Julia set was proved in [2,4,7]. Computer generated images are given in [2].\\

     THEOREM 2. If $f_c(z) = z^2 + c$ with hyperbolic numbers $z = x + jy, c = a + jb$, then
        $\mathbf{K}_{\mathbf{H}}(f_c)$ has one the following four possibilities depending on the characteristic coordinates  $c_1 = a - b , c_2 = a + b$ of $c$:\\      
     
      (i) $\mathbf{K}_{\mathbf{H}}(f_c)$ is connected if $c\in S$, i.e. if $c_1, c_2 \in [-2, \frac{1}{4}]$, namely it is a rectangle. \\    
   
   If $c\notin S$:\\
      
      (ii) $\mathbf{K}_{\mathbf{H}}(f_c)$ is a Cantor set, thus totally disconnected, if $c_1, c_2 < -2$\\
      
      (iii) $\mathbf{K}_{\mathbf{H}}(f_c)$ is disconnected but not totally disconnected
      if either $c_1$ or $c_2$ is in $[-2, \frac{1}{4}]$ and the other is $<-2$\\
     
      (iv)  $\mathbf{K}_{\mathbf{H}}(f_c)$ is empty otherwise, i.e. if  either $c_1$ or $ c_2  > \frac{1}{4}$ \\
           
   PROOF: (i) Using Lemma 1 and Lemma 2, if $c\in S$, then  $\mathbf{K}_{\mathbf{H}}(f_c) $ is an ordered pair of the two intervals  $ [ -p_+^k, p_+^k] $  where      
 $p_+ ^k $ is the positive fixed point of the real quadratic function $  x^2 + c_k, k = 1,2$ .\\
      (ii) An ordered pair of two Cantor sets is a Cantor set.\\
      (iii) An ordered pair of a connected set and a totally disconnected set is a disconnected set.\\
      (iv) At least one of the real filled Julia sets  $\mathbf{K}_{\mathbf{R}} (f_{c_k}) $ is empty for $k=1,2$, and thus the pair.
       
     \vspace{0.3cm}  
     
     OPEN QUESTIONS\\
     (1) What are the properties of other maps, which are important in complex analysis, if they are considered as maps over hyperbolic numbers? One such map has been treated in [8], namely the generalization of the complex Riemann zeta function $\zeta(s) =  \sum_{n=1}^{\infty}\frac{1}{n^s},$ s a complex number, to hyperbolic numbers $s\in \mathbf{H}$. In [8] it is shown that every zero of the hyperbolic Riemann zeta function is trivial, so there is no analogy to the Riemann hypothesis for hyperbolic numbers. \\
     (2) What maps over the hyperbolic numbers have interesting dynamics, if any?
       
    \vspace{0.3cm}
    REFERENCES\\
   
    [1]  Artzy, R., 1992, Dynamics of quadratic functions in cycle planes. $\Large {Journal\ of\ Geometry}$ 44, 26-32.\\
   
    [2]  Blankers, V.,Rendfrey,T., Shukert, A. and Shipman, P.D.,  2019, Julia and Mandelbrot Sets for Dynamics over the Hyperbolic Numbers. $\Large{Fractal \ Frac}$        3,6.\\
    
    [3]  Devaney, R., 1996, A First Source in Chaotic Dynamical Systems: Theory and Experiment, Westview Press.\\
    
    [4]  Fishback, P.E.,  2005, Quadratic dynamics in binary number systems. $\Large {Journal\ of \ Difference \ Equations \ and \ Applications}$,Vol.11, No. 7, 597-603.\\
   
    [5]  Metzler, W., 1994, The "Mystery" of the quadratic Mandelbrot set. $\Large{American \ Journal\ of \ Physics }$ 62, 813-814.\\
   
    [6] Rochon, D., 2000, A generalized Mandelbrot set for bicomplex numbers. $\Large{Fractals}$ 8, 355-368. \\
   
    [7] Rochon, D., 2003, On A Generalized Fatou-Julia Theorem. $\Large{Fractals}$ 11, 3, 213-219. \\
   
    [8] Rochon, D., 2004, A  Bicomplex Riemann Zeta Function. $\Large{Tokyo\  Journal \ Math}$, 27, 3, 357-369.\\
   
    [9] Cantoni, F., Boccaletti, D., Cannata, R., Catoni, V., Nichelatti, E. and Zampetti, P., 2008, The Mathematics of Minkowski Space-Time. Birkhaeuser Verlag.\\

         Department of Mathematics,The Graduate Center of the City University of New York,   shayes@gc.cuny.edu.  \\

      \end{document}